\documentclass[graybox]{svmult}

\usepackage{mathptmx}        
\usepackage{makeidx}         
\usepackage{cite}
\usepackage{graphicx}        
\usepackage{grffile}
\usepackage{epstopdf}        
\usepackage{amsmath}
\usepackage{multicol}        
\usepackage[bottom]{footmisc}

\makeindex             

\begin{document}


\title*{Hundreds of new satellites of figure-eight orbit  computed with  high precision}

\author{I. Hristov$^{1,a}$, R. Hristova$^1$, I. Puzynin$^2$,  T. Puzynina$^2$, Z. Sharipov$^{2,b}$, Z. Tukhliev$^2$\\
\vspace{0.5cm}
\emph{$^1$ Faculty of Mathematics and Informatics, Sofia University, Sofia, Bulgaria}\\
\emph{$^2$ Meshcheryakov Laboratory of Information Technologies, Joint Institute for Nuclear Research,  Dubna, Russia}\\
\vspace{0.5cm}
\emph{E-mails:  $^a$  ivanh@fmi.uni-sofia.bg \hspace{0.3cm} $^b$ zarif@jinr.ru}}
\authorrunning{I. Hristov, R. Hristova, I. Puzynin et al.}
\maketitle


\abstract {Satellites (topological powers) of the famous figure-eight orbit are special periodic solutions of the planar three-body problem.
           In this paper we use a modified Newton's method based on the Continuous analog of  Newton's method and high precision arithmetic
           for a purposeful numerical search of new satellites of the figure-eight orbit.
           Over 700 new satellites are found, including 76 new linearly stable ones. 7 of the newly found linearly stable satellites are choreographies.
           The linear stability is checked by a high precision computing of the eigenvalues of the monodromy matrices.
           The initial conditions of all found solutions are given with 150 correct decimal digits.}

\section{Introduction}

Numerical algorithms used in finding new periodic orbits are of big importance
for further investigation and understanding of the three-body problem.
Major progress in this area has been made in recent years.
In 2013 Shuvakov and Dmitrashinovich found 13 new topological families of periodic solutions for the planar three-body problem,
applying a clever numerical algorithm based on the gradient descent method in the standard double precision arithmetic \cite{Shuv1, Shuv2}.
Since the three-body problem has a sensitive dependence
on the initial conditions,  working with double precision limits the number of solutions that can be found,
especially nonstable solutions. Nevertheless, a method using double precision arithmetic can still be very useful, since it is not time consuming and catches many of the important stable solutions.

In 2017 Li and Liao applied the Newton's method to search for new periodic orbits of the planar three-body problem.
They used the high order  multiple precision Taylor series method to form the linear system at each step of Newton's method
and found  more than 600 new topological families of periodic orbits \cite{Liao1}. This result shows the great potential of high precision
computations for the three-body problem.
The absence of details for the used numerical algorithm
in \cite{Liao1} motivated us to present this algorithm in our recent work  \cite{Ourpaper}. In \cite{Ourpaper} we made a detailed description of the Newton's method and also its modification, based on the Continuous analog of Newton's method \cite{CANM} for computing periodic orbits of the planar three-body problem. In particular, we gave all needed formulas, based on the rules of the automatic differentiation for the Taylor series method. We tested our programs by making a general search for relatively short periods and a relatively coarse search grid. As a result we found 33 new topological families, that are not included in the database in \cite{Liao1}. The solutions were computed with 150 correct digits. We concluded that the usage of a relaxed return proximity condition for the initial approximations and the usage of the modified Newton's method instead of the classic Newton's method is a very perspective approach to find more periodic solutions for a given search grid.

In this paper we use the presented in \cite{Ourpaper} modified Newton's method for a purposeful search of new satellites (topological powers) of the famous figure-eight orbit \cite{Moore}.
There are two known figure-eight solutions: Moore's figure eight orbit \cite{Moore, Montgomery}, which is a stable choreography,
and Simo's figure-eight orbit which is of hyperbolic-elliptic type and is not a choreography \cite{Simo}.
A choreography is such a periodic solution, for which all bodies move along one and the same trajectory with a time delay of $T/3$, where $T$ is the period of the solution.
Satellites of figure-eight are of special interest as it is expected that some of them are stable periodic orbits and (or) choreographies.
Our work is motivated by the works \cite{Shuv3, Shuv4}, where many new satellites of figure-eight with zero angular momentum are found, including new choreographies,
and particularly one stable choreography, which is until now only the second stable choreography found after the Moore's figure-eight orbit.
Many three-body choreographies (345) are also found in \cite{Simo}, but they are with nonzero angular momentum and undetermined topological type.

In this work we first make a general search and then concentrate on a small domain for initial conditions. As a result we manage
to find 821 satellites of the figure-eight orbit (778 new ones), including 7 new linearly stable choreographies.

\section{Mathematical model}
A planar motion of the three bodies is considered. The bodies are with equal masses and are treated as mass points.
The normalized differential equations describing the motion of the bodies are:
\begin{equation}
\ddot{r}_i=\!\!\!\sum_{j=1,j\neq i}^{3} \frac{(r_j-r_i)}{{\|r_i -r_j\|}^3}, \quad i=1,2,3.
\end{equation}
The vectors $r_i$, $\dot{r}_i$ have two components: $r_i=(x_i, y_i)$, $\dot{r}_i=(\dot{x}_i, \dot{y}_i)$.
The system (1) can be written as a first order one this way:
\begin{equation}
\dot{x}_i={vx}_i, \hspace{0.15 cm} \dot{y}_i={vy}_i, \hspace{0.15 cm} \dot{vx}_i=\!\!\!\sum_{j=1,j\neq i}^{3} \frac{(x_j-x_i)}{{\|r_i -r_j\|}^3}, \hspace{0.15 cm} \dot{vy}_i=\!\!\!\sum_{j=1,j\neq i}^{3} \frac{(y_j-y_i)}{{\|r_i -r_j\|}^3}, \quad i=1,2,3
\end{equation}
We solve numerically the problem in this first order form. Hence we have a vector of 12 unknown functions
$ X(t)={(x_1, y_1, {vx}_1, {vy}_1, x_2, y_2, {vx}_2, {vy}_2, x_3, y_3, {vx}_3, {vy}_3)}^\top$.
We search for periodic planar collisionless orbits as in \cite{Shuv1, Shuv2, Liao1}: with zero angular momentum and symmetric initial
configuration with parallel velocities:
\begin{equation}
\begin{aligned}
(x_1(0),y_1(0))=(-1,0), \hspace{0.2 cm} (x_2(0),y_2(0))=(1,0), \hspace{0.2 cm} (x_3(0),y_3(0))=(0,0) \\
({vx}_1(0),{vy}_1(0))=({vx}_2(0),{vy}_2(0))=(v_x,v_y) \hspace{2 cm}\\
({vx}_3(0),{vy}_3(0))=-2({vx}_1(0),{vy}_1(0))=(-2v_x, -2v_y) \hspace{1.5 cm}
\end{aligned}
\end{equation}
Here the velocities $v_x\in [0,1], v_y\in [0,1]$ are parameters. The periods of the orbits are denoted with $T$. So, our goal is to find triplets $(v_x, v_y, T)$
for which the periodicity condition $X(T)=X(0)$ is fulfilled.

\section{Description of the numerical searching procedure}
The numerical searching procedure for new periodic orbits has three stages. During the first stage we search for candidates for correction
with the modified Newton's method, i.e. we compute initial approximations.
During the second stage we apply the modified Newton's method, which has a larger domain of convergence
than the classic Newton's method.
Convergence during the second stage means that a periodic orbit is found.
During the third stage  we apply the classic Newton's method with a higher precision
in order to specify the solutions with more correct digits (150 correct digits in this work).

A square 2D grid with some step is introduced in some 2D domain for parameters
$v_x, v_y$ for the first stage. For every grid point the system (2) with initial conditions (3) is solved numerically up to a pre-fixed time
$T_0$. Let $$P(t)=\sqrt{\sum_{i=1}^{3}{\|r_i(t)-r_i(0)\|}^{2}+\sum_{i=1}^{3}{\|\dot{r}_i(t)-\dot{r}_i(0)\|}^{2}}, \hspace{0.5 cm} 0 \leq t \leq T_0$$ be the proximity function.
As candidates for the correction method, we take those triplets  $(v_x, v_y, T)$ for which the return proximity $P(T)$ is less then $0.1$:
$$P(T) = \min_{1<t\leq T_0}P(t) < 0.1$$
We also set the constraint that the return proximity has a local minimum on the grid for $v_x, v_y$.
The modified and the classic Newton's method are applied for correction during the second and the third stage.
To do this, a linear algebraic system with a $12\times3$ matrix for the corrections $\Delta v_x,
\Delta v_y, \Delta T$ has to be solved \cite{Abad}. To form the matrix, a system of 36 ODEs (the original differential equations plus the differential equations for the partial derivatives with respect
to the parameters $v_x$ and $v_y$) has to be solved first.  The linear algebraic system is solved in linear least square sense
using $QR$ decomposition based on Householder reflections \cite{Demmel}. The multiple precision Taylor series method \cite{Barrio1, Barrio2}
with a variable step-size strategy from \cite{Jorba} is used during all three stages to solve the corresponding systems of ODEs.
During the first and the second stage Taylor series method with order 154 and 134 decimal digits of precision is used.
During the third stage we make two computations. The first computation is with 220-th order Taylor series method and 192
digits of precision and the second computation is for verification - with 264-th order method and 231 digits of precision.
All the details of the Newton's and the modified Newton's method for
the planar three-body problem can be seen in our recent work  \cite{Ourpaper}.

A purposeful search for figure-eight satellites is made in this work.
We first make a general search in the rectangle $[0, 0.8]\times [0, 0.8]$ for $v_x, v_y$
with a grid step $1/2048$ and $T_0=70$. After that we add to our found satellites the satellites found in \cite{Shuv3, Shuv4, Dimi}
(we use the tables in \cite{Shuv3, Shuv4} and the table for sequence V (figure-eight) at \cite{Table}).
A cluster of the initial velocities' points $(v_x, v_y)$ is seen, which we enclose with the rectangle $[0.1, 0.33]\times [0.49, 0.545]$ (see Fig.~1).
All the points outside of this rectangle are surrounded by circles with radius $0.02$.
Then  we concentrate on this rectangle and the union of circles with a grid step $1/4096$ and $T_0=300$.
After that we left only those points outside the rectangle $[0.1, 0.33]\times [0.49, 0.545]$ for which
we found in their vicinities (circles) some new solutions. Then we again consider the union of vicinities  of the left points together with the newly found solutions,
i.e. we make a second iteration for searching. All in all we make a search for the rectangle $[0.1, 0.33]\times [0.49, 0.545]$
plus the domain with the curved boundary shown in Fig.~1.

\begin{figure}
\begin{center}
\includegraphics[scale=0.58]{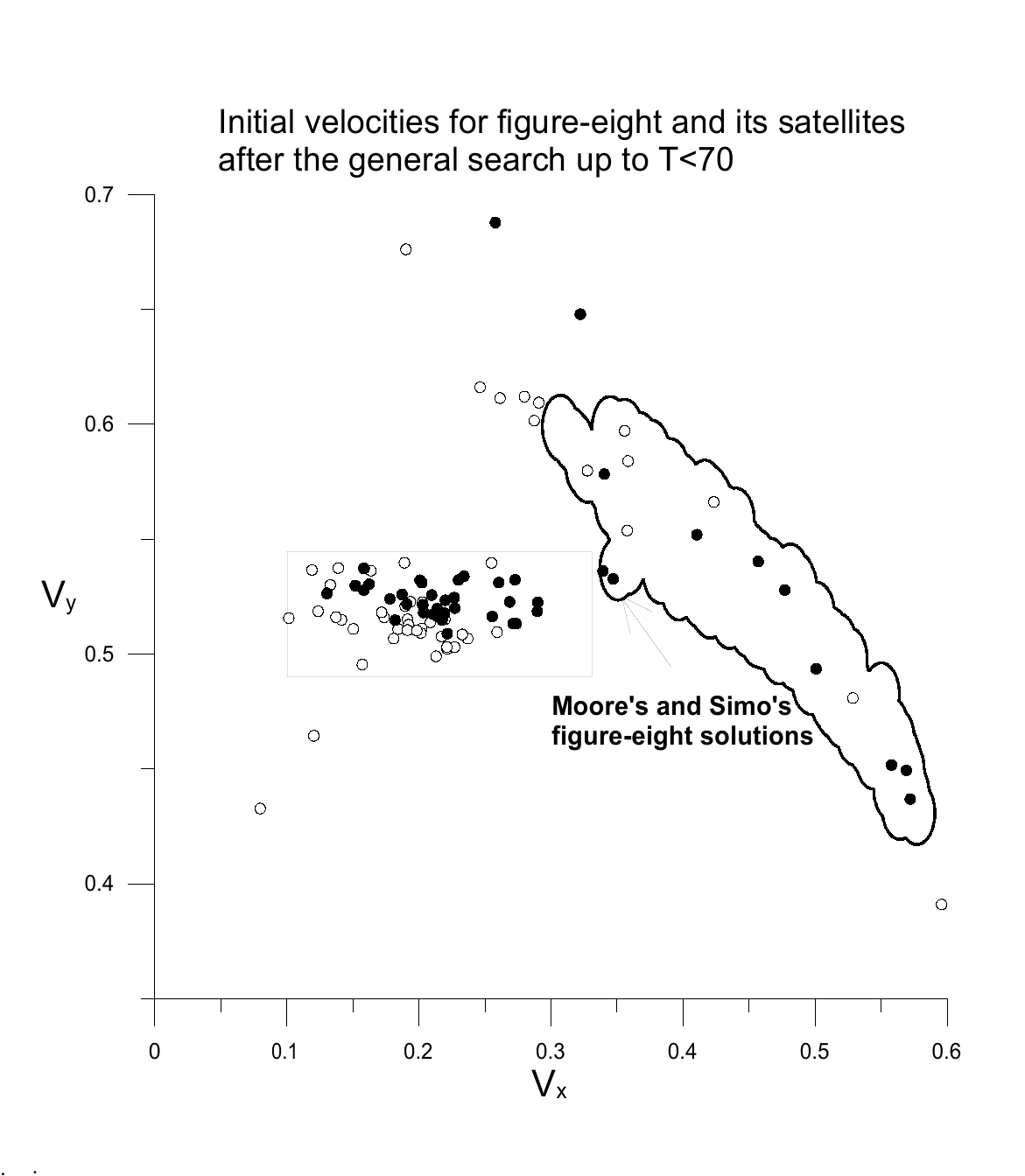}
\caption{Initial velocities: 45 black points - previously found satellites, 50 white points - newly found satellites after a general search up to $T_0=70$ with a grid-step 1/2048}
\label{fig:1}
\end{center}
\end{figure}

\section{Linear stability computations}
The linear stability of a given periodic orbit $X(t), 0 \leq t \leq T,$ is determined by the eigenvalues $\lambda$ of the $12 \times 12 $ monodromy matrix $M(X;T)$ \cite{Roberts}:
$$M_{ij}[X;T] = \frac{\partial X_i(T)}{\partial X_j(0)}, \quad  M(0)=I $$
The eigenvalues of $M$ come in pairs or quadruplets: $(\lambda, \lambda^{-1}, \lambda^{*}, \lambda^{*-1})$. They are of four types:

1) Elliptically stable - $\lambda = \exp(\pm 2 \pi i \nu),$ where $\nu > 0$ (real) is the stability angle. In this case the eigenvalues are on the unit circle.
Angle $\nu$ describes the stable revolution of adjacent trajectories around a periodic orbit.

2) Marginally stable - $\lambda = \pm 1.$

3) Hyperbolic - $\lambda = \pm\exp(\pm \mu),$ where $\mu>0$ (real) is the Lyapunov exponent.

4) Loxodromic - $\lambda = \exp(\pm\mu \pm i \nu),$  $\mu, \nu$ (real)

Eight of the eigenvalues of $M$ are equal to $1$ \cite{Roberts}. The other four determine the linear stability.
Here we are interested in elliptically stable orbits, i.e. the four eigenvalues to be $\lambda_j= \exp(\pm 2 \pi i \nu_j), j=1,2$.
The elements of the monodromy matrix $M_{ij}$ are computed in the same way as the partial derivatives in \cite{Ourpaper}
- with the multiple precision Taylor series method using the rules of automatic differentiation.
For computing the eigenvalues we use a Multiprecision Computing Toolbox \cite{Advanpix} for MATLAB \cite{Matlab}.
First the elements of $M$ are obtained with 120 correct digits and then two computations  with 70 and 120 digits of precision are made
with the toolbox. The four eigenvalues under consideration  are verified by a check for matching the first 30 digits of them and the corresponding
condition numbers obtained by the two computations (with 70 and 120  digits of precision).
The obtained values for Moore's figure-8 are {\footnotesize $\nu_1=0.298092529004750122423759217553$}, {\footnotesize $\nu_2=0.00842272470813137798255957006197$}
and for hyperbolic-elliptic Simo's figure-8 are {\footnotesize $\lambda=1.09587441937857035881716297118$}, {\footnotesize $\nu=0.297572835261762177616886515603$}.
First digits exactly match with given in \cite{Simo}.

\section{Results}
We found totally 821 figure-eight satellites
(including Moore's and Simo's figure-eight orbits and old satellites, and counting different initial conditions as different solutions).
778 of the 821 solutions are new (not included in \cite{Shuv3, Shuv4, Dimi}). 65 of them are choreographies (45 new choreographies).

A topological method from \cite{Mont} is applied to classify the periodic orbits into families.
Each family corresponds to a different conjugacy  class of the free group on two letters $(a,b)$.
We use "the free group word reading algorithm" from \cite{Shuv2} to obtain the free group elements. Satellites of figure-eight correspond to free group elements ${(abAB)}^k$
for some natural power $k$. For the found solutions the power $k$ varies from 4 up to 62, including all values in between.
For each found solution we computed the free group element and the four numbers $(v_x,v_y,T,T^*)$ with 150 correct digits,
where $T^*$ is the scale-invariant period. The scale-invariant period is defined as  $T^{*}=T{|E|}^{\frac{3}{2}},$
where  $E$ is the energy of our initial configuration: $E=-2.5 + 3({v_x}^2+{v_y}^2)$. Equal $T^*$ for two different initial conditions means two different representations of the same solution.
All 821 solutions  computed with 150 correct digits, plots of the trajectories and animations in the real $x-y$ plane for them can be seen in \cite{rada3body}.

We started the investigation for the stability of the presented solutions with simple tests.
The simple tests consist of integration of system (2) for the found initial conditions with Taylor series method with order 70 and 57 decimal digits of precision.
The solutions that do not escape from the periodic trajectory for a very long time (2000 periods)
become candidates for stable solutions. We obtain 115 candidates. Let us mention that if the greatest magnitude of eigenvalues is  greater  but very closed to one,
these tests may give false result, i.e. giving that a solution is stable, while it is not.
By computing the eigenvalues of the monodromy matrices we obtain that actually we have 82 linearly stable orbits (76 new ones).
The obtained linearly stable solutions are only among the 115 candidates, which passed the first simple test.
The linear stability results for the earlier found satellites  from \cite{Dimi} are confirmed, i.e.
we obtain the same linear stability results for old solutions as in \cite{Dimi}.

The results for linear stability can be seen in a table in \cite{rada3body}.
All stability angles $\nu_{1,2}$ with 30 correct digits for the linearly stable solutions
are given in another table in \cite{rada3body}.
The distribution of the initial condition points
can be seen in Fig.~2 (the black points are the linearly stable solutions, the white points -- the unstable ones (more precisely, not confirmed to be linearly stable)).
As it is seen from this figure, most of the initial condition points for the linearly stable solutions are found in the rectangle $[0.1, 0.33]\times [0.49, 0.545]$
and also these points are clustered.

\begin{figure}
\begin{center}
\includegraphics[scale=0.7]{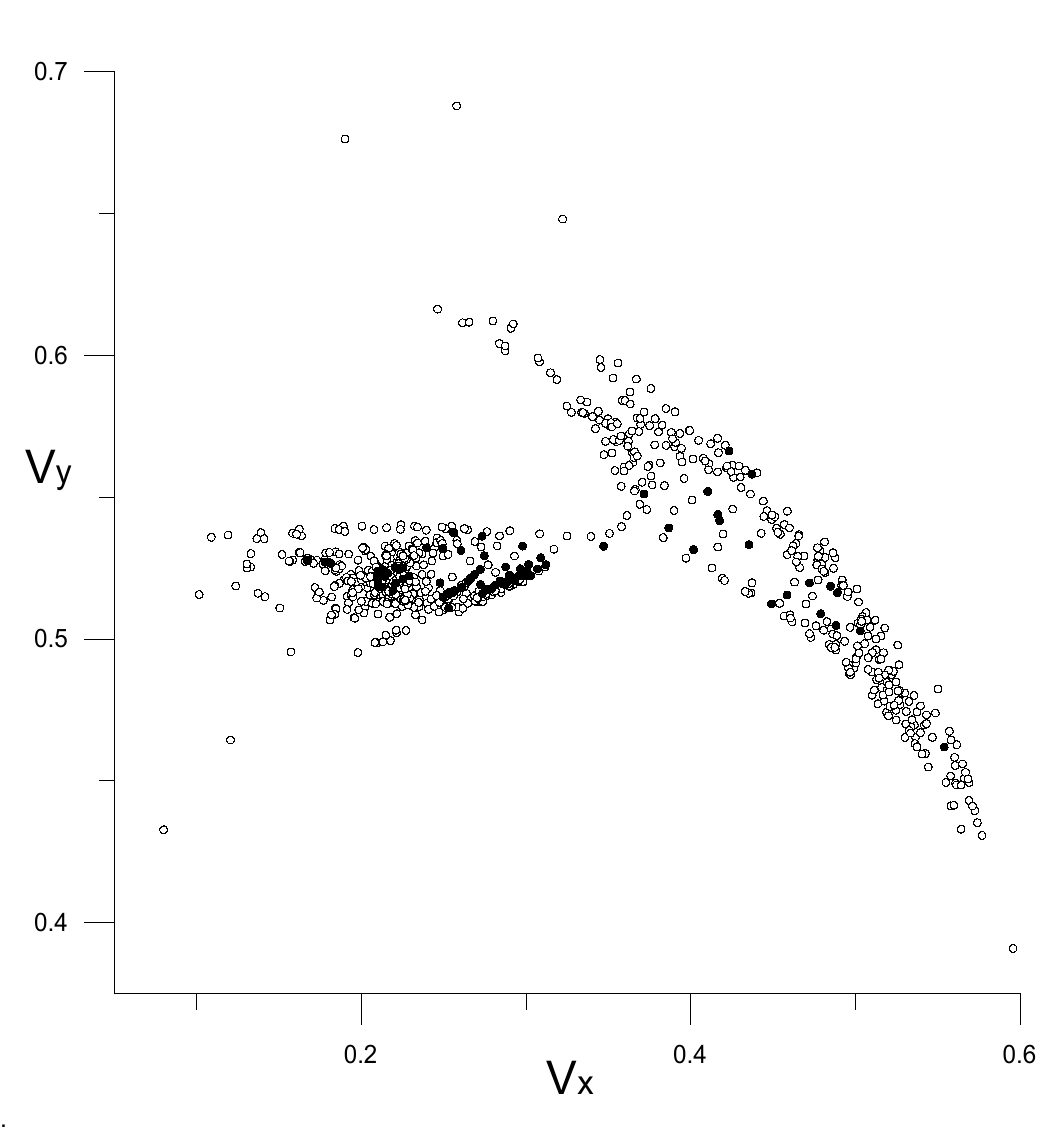}
\caption{Initial velocities for all found solutions (including previously found), black points are the linearly stable solutions,
white points - unstable (more precisely, not confirmed to be linearly stable) }
\label{fig:1}
\end{center}
\end{figure}

 Moore's and Simo's figure-eight form a pair of solutions with equal topological power and very close $T^{*}$ where one orbit is linearly stable, and the other is of hyperbolic-elliptic type.
 The two solutions have the following property: "The stability angle of the elliptic eigenvalues of the hyperbolic-elliptic type solution is very close to the larger stability angle
 of the linearly stable solution and the larger hyperbolic eigenvalue of the hyperbolic-elliptic type solution is greater than one but very close to one".
 11 such pairs are observed in the found set of solutions. For example, there exists a pair of solutions with {\footnotesize $\nu_1=0.275011350174751332178664212923$}, {\footnotesize $\nu_2=0.00341516982815402730209182950586$} and {\footnotesize $\nu=0.274944611347846047797427400637$}, {\footnotesize $\lambda=1.02168877842035886643036103128$}. The initial conditions for this pair can
 be seen in Table 1. The trajectories of the three bodies in the real $x-y$ plane can be seen in Fig.~3 and Fig.~4. A table with the data for all the 11 pairs can be seen in
 \cite{rada3body}. The only pair for which one of the solution is choreography and the other is not, is Moore's and Simo's figure-8 pair.
 The other pairs include solutions which are both choreographies or both are not.

 \begin{table}
 \begin{tabular}{ p{0.3cm} p{2.75cm} p{2.75cm} p{2.75cm} p{2.75cm} p{0.3cm} }
 \hline
 $N$ & $\hspace{1.1cm}v_x$ & $\hspace{1.1cm}v_y$ & $\hspace{1.2cm}T$ & $\hspace{1.2cm}T^*$ & $\hspace{0.05cm}k$ \\
 \hline
\scriptsize{1} & \scriptsize{0.43704391669569526e0} & \scriptsize{0.55809267787922151e0} & \scriptsize{0.25199131683991163e3} & \scriptsize{0.24919018052039563e3} & \scriptsize{27} \\
\scriptsize{2} & \scriptsize{0.54316996361672422e0} & \scriptsize{0.4733330189162051e0} & \scriptsize{0.27222283448647228e3}  & \scriptsize{0.24919018052115745e3} & \scriptsize{27} \\
\end{tabular}

{\caption {Data with 17 correct digits for a pair of linearly stable and hyperbolic-elliptic solutions}}
\end{table}

 \begin{figure}
\begin{center}
\includegraphics[scale=0.52]{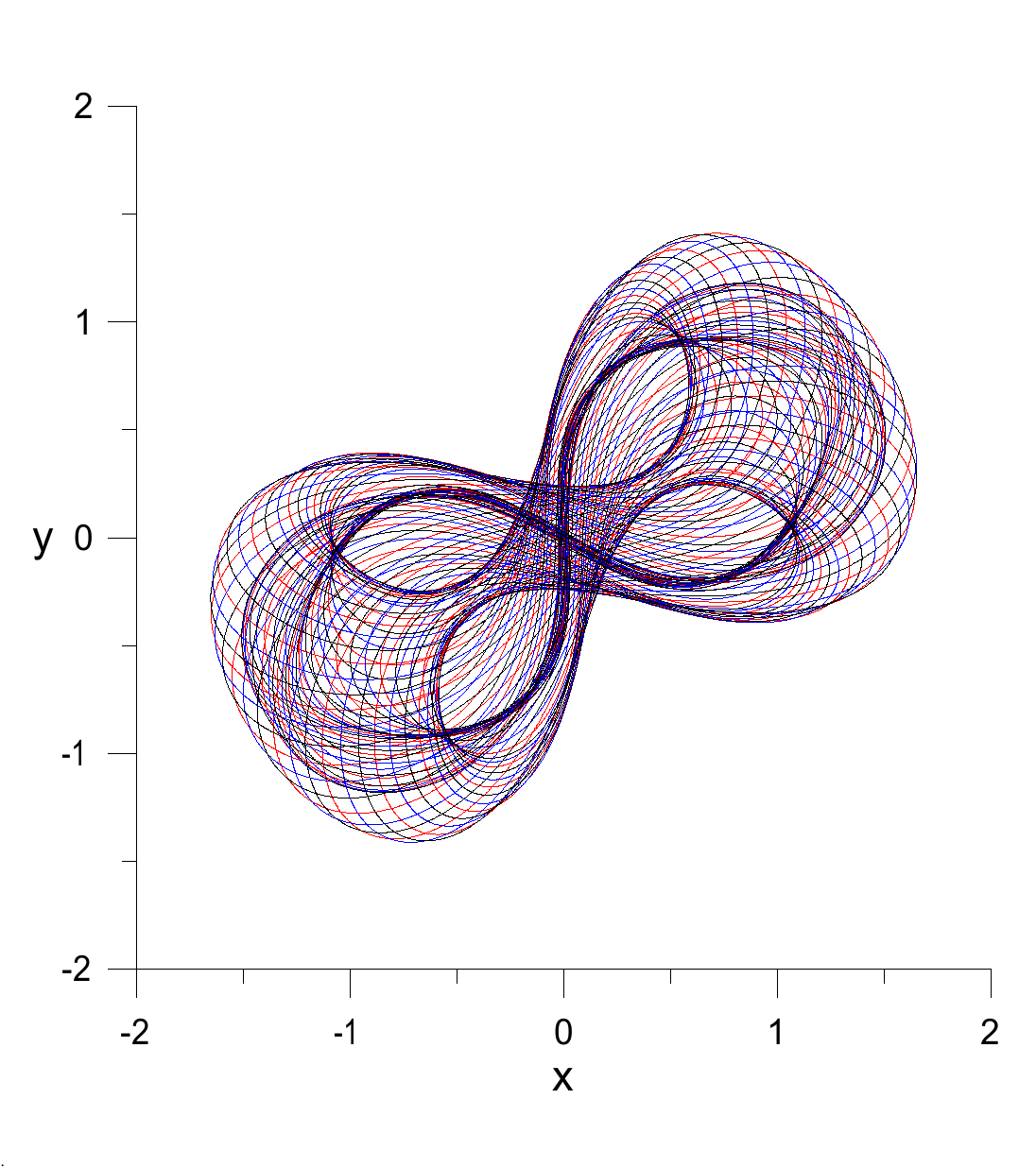}
\caption{Linearly stable solution (solution 1 from Table 1)}
\label{fig:1}
\end{center}
\end{figure}

\begin{figure}
\begin{center}
\includegraphics[scale=0.52]{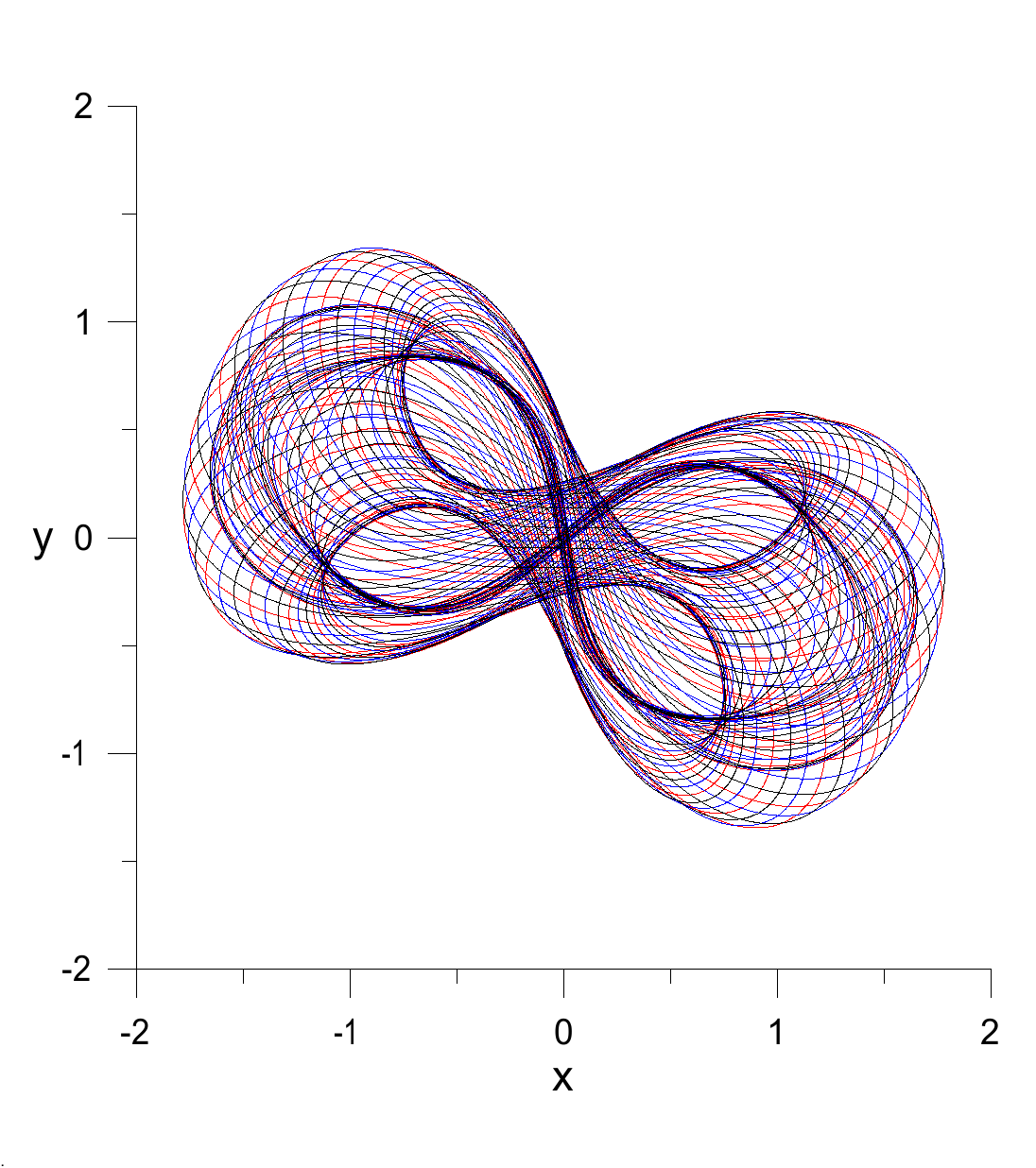}
\caption{Hyperbolic-elliptic type solution (solution 2 from Table 1)}
\label{fig:1}
\end{center}
\end{figure}

An important result is the finding of 7 new linearly stable choreographies which are topological powers of the figure-eight.
The second and the third of the choreographies are presented by two different initial conditions (all in all we have 9 solutions).
The numbers $(v_x,v_y,T,T^*)$ with 17 correct digits and the power $k$ can be seen in Table 2.
The stability angles are given in Table 3.
Until now there were known only two stable choreographies with zero angular momentum - the famous Moore's figure-eight and the Shuvakov's choreography \cite{Shuv3}
(solutions 28 and 29 in the table in \cite{rada3body}).
The first linearly stable choreography form a pair (as explained earlier) with a hyperbolic-elliptic choreography with
{\footnotesize $\lambda= 1.00000000099138284038597245085$}, {\footnotesize $\nu=0.48697042062642601225914847037$}. The trajectory in the real $x-y$ plane of the hyperbolic-elliptic choreography corresponding to the first linearly stable choreography
can be seen in Fig.~5. The trajectories of the linearly stable choreographies can be seen in Figures 6 -- 12.

The extensive computations for the numerical search are performed in "Nestum" cluster, Sofia, Bulgaria \cite{nestum} and "Govorun" supercomputer, Dubna, Russia \cite{HybriLIT}.
The GMP library (GNU multiple precision library) \cite{GMP} for multiple precision floating point arithmetic is used.

 \begin{table}
 \begin{tabular}{ p{0.3cm} p{2.75cm} p{2.75cm} p{2.75cm} p{2.75cm} p{0.3cm} }
 \hline
 $N$ & $\hspace{1.1cm}v_x$ & $\hspace{1.1cm}v_y$ & $\hspace{1.2cm}T$ & $\hspace{1.2cm}T^*$ & $\hspace{0.05cm}k$ \\
 \hline
\scriptsize{1} & \scriptsize{0.40169296201102497e0} & \scriptsize{0.53156253772635003e0} & \scriptsize{0.2706819009639754e3} & \scriptsize{0.34179371490016304e3} & \scriptsize{37} \\
\scriptsize{2} & \scriptsize{0.29791817208757154e0} & \scriptsize{0.53275520940835742e0} & \scriptsize{0.21032180603055653e3} & \scriptsize{0.34179379656326958e3} & \scriptsize{37} \\
\scriptsize{3} & \scriptsize{0.30889168505067089e0} & \scriptsize{0.52866005593876979e0} & \scriptsize{0.2119149795956615e3} & \scriptsize{0.34179379656326958e3} & \scriptsize{37} \\
\scriptsize{4} & \scriptsize{0.17791680921792612e0} & \scriptsize{0.52719348724657283e0} & \scriptsize{0.17814022614781591e3} & \scriptsize{0.35085277996057519e3} & \scriptsize{38} \\
\scriptsize{5} & \scriptsize{0.2529928848686255e0} & \scriptsize{0.51092018996687722e0} & \scriptsize{0.18632778573173244e3} & \scriptsize{0.35085277996057519e3} & \scriptsize{38} \\
\scriptsize{6} & \scriptsize{0.27224186760209537e0} & \scriptsize{0.52452018883171169e0} & \scriptsize{0.23223664454023032e3} & \scriptsize{0.40645290803177171e3} & \scriptsize{44} \\
\scriptsize{7} & \scriptsize{0.29653206621705968e0} & \scriptsize{0.52476506719403896e0} & \scriptsize{0.24274490570242023e3} & \scriptsize{0.40645421334926255e3} & \scriptsize{44} \\
\scriptsize{8} & \scriptsize{0.21036388572856618e0} & \scriptsize{0.51854057011033809e0} & \scriptsize{0.27480620837192349e3} & \scriptsize{0.53574606807894121e3} & \scriptsize{58} \\
\scriptsize{9} & \scriptsize{0.26684473313570582e0} & \scriptsize{0.52168674682029369e0} & \scriptsize{0.30063538172570104e3} & \scriptsize{0.53576756803804106e3} & \scriptsize{58} \\

\end{tabular}
{\caption {Data with 17 correct digits for the 7 new linearly stable choreographies}}
\end{table}

\vspace{- 1 cm}

\begin{table}
 \begin{tabular}{ p{1 cm} p{5.5 cm} p{5.5 cm}}
 \hline
 $N$ & $\hspace{2 cm}\nu_1$ & $\hspace{2 cm} \nu_2$\\
 \hline
  \scriptsize{1} & \scriptsize{0.486970420626426013258486704454} & \scriptsize{1.57783479465699085471836250548e-10}\\
  \scriptsize{2} & \scriptsize{0.0327953388952792073716032366859} & \scriptsize{3.22317484874566580526492303709e-06}\\
  \scriptsize{3} & \scriptsize{0.0327953388952792073716032366859} & \scriptsize{3.22317484874566580526492303709e-06}\\
  \scriptsize{4} & \scriptsize{0.267131864119191530384832354753} & \scriptsize{0.127489871368591999957771257333}\\
  \scriptsize{5} & \scriptsize{0.267131864119191530384832354753} & \scriptsize{0.127489871368591999957771257333}\\
  \scriptsize{6} & \scriptsize{0.187691189591235263346007262414} & \scriptsize{2.61504980535709161857184424983e-12}\\
  \scriptsize{7} & \scriptsize{0.0937385311958221201020392196211} & \scriptsize{0.000838678707197693244105921997421}\\
  \scriptsize{8} & \scriptsize{0.427373702369188544983378448257} & \scriptsize{0.00236826906538802133452909107108}\\
  \scriptsize{9} & \scriptsize{0.376315052261497474047661511684} & \scriptsize{4.95737788455128082115378851449e-10}\\

\end{tabular}
{\caption {Stability angles $\nu_j$, $\lambda_j = \exp(\pm2\pi i \nu_j)$, j=1,2 for the 7 new linearly stable choreographies}}
\end{table}

\section{Conclusions}
A  modified Newton's method, based on the Continuous analog of  Newton's method and high precision is successfully used in finding more than 700
new satellites of the famous figure-eight orbit, including 76 new linearly stable ones.  7 of the newly found linearly stable satellites are choreographies.
The results show the great potential of the used numerical algorithm in
searching for new periodic orbits of the three-body problem.

\begin{acknowledgement}

We thank for the opportunity to use the computational resources of the "Nestum" cluster, Sofia, Bulgaria
and  the "Govorun" supercomputer at the Meshcheryakov Laboratory of Information Technologies of JINR, Dubna, Russia.
We would also like to thank Veljko Dmitrashinovich from Institute of Physics, Belgrade University, Serbia for a valuable e-mail discussion and advice,
and his encouragement to continue our numerical search for new periodic orbits.
The work is supported by a grant of the Plenipotentiary Representative of the Republic of Bulgaria at JINR, Dubna, Russia.

\end{acknowledgement}

\begin{figure}
\begin{center}
\includegraphics[scale=0.52]{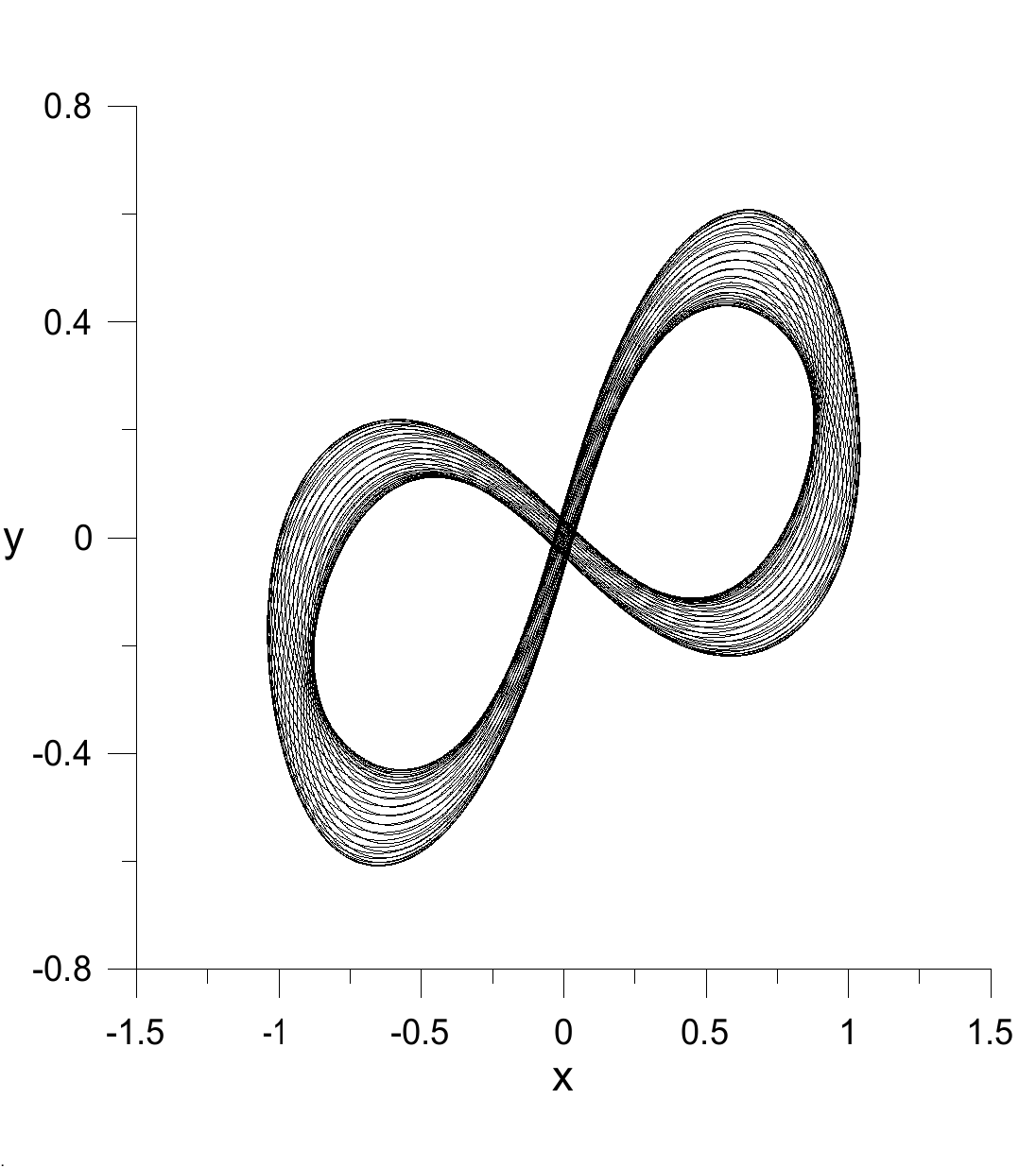}
\caption{Hyperbolic-elliptic  choreography corresponding to the stable choreography (solution 1)}
\label{fig:1}
\includegraphics[scale=0.52]{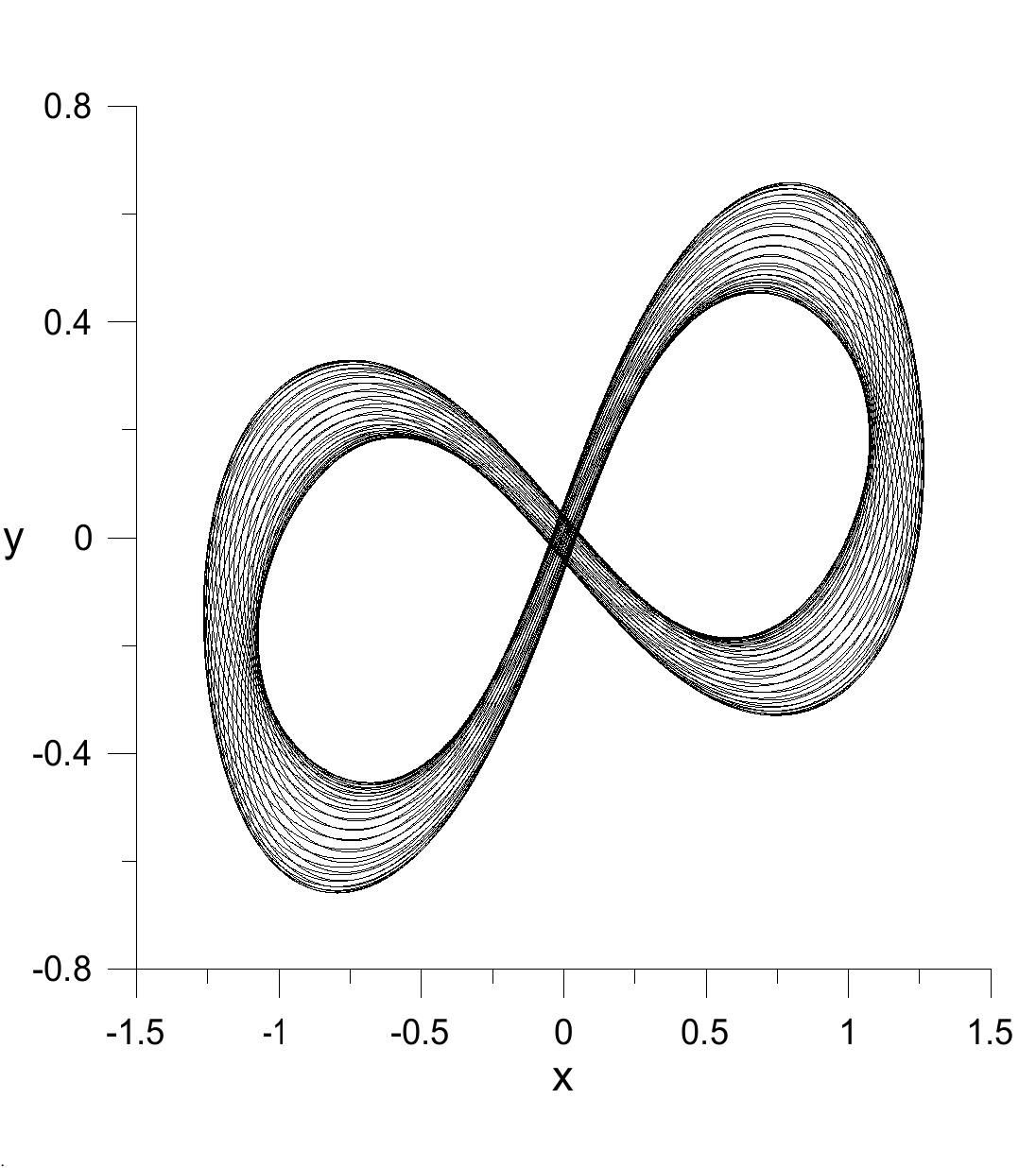}
\caption{Linearly stable choreography (solution 1)}
\label{fig:1}
\end{center}
\end{figure}

\begin{figure}
\begin{center}
\includegraphics[scale=0.52]{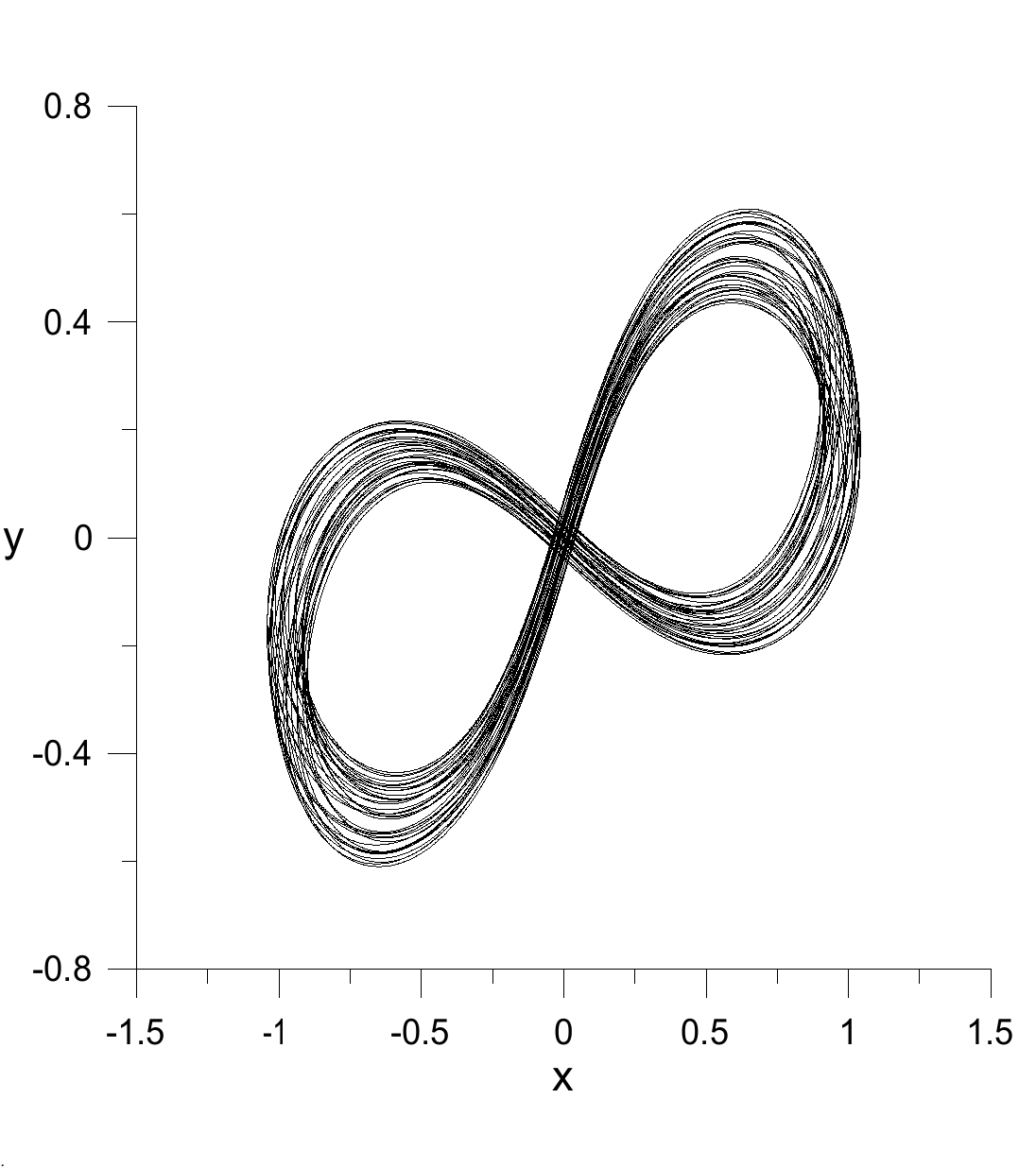}
\caption{Linearly stable choreography (solution 2)}
\label{fig:1}
\includegraphics[scale=0.52]{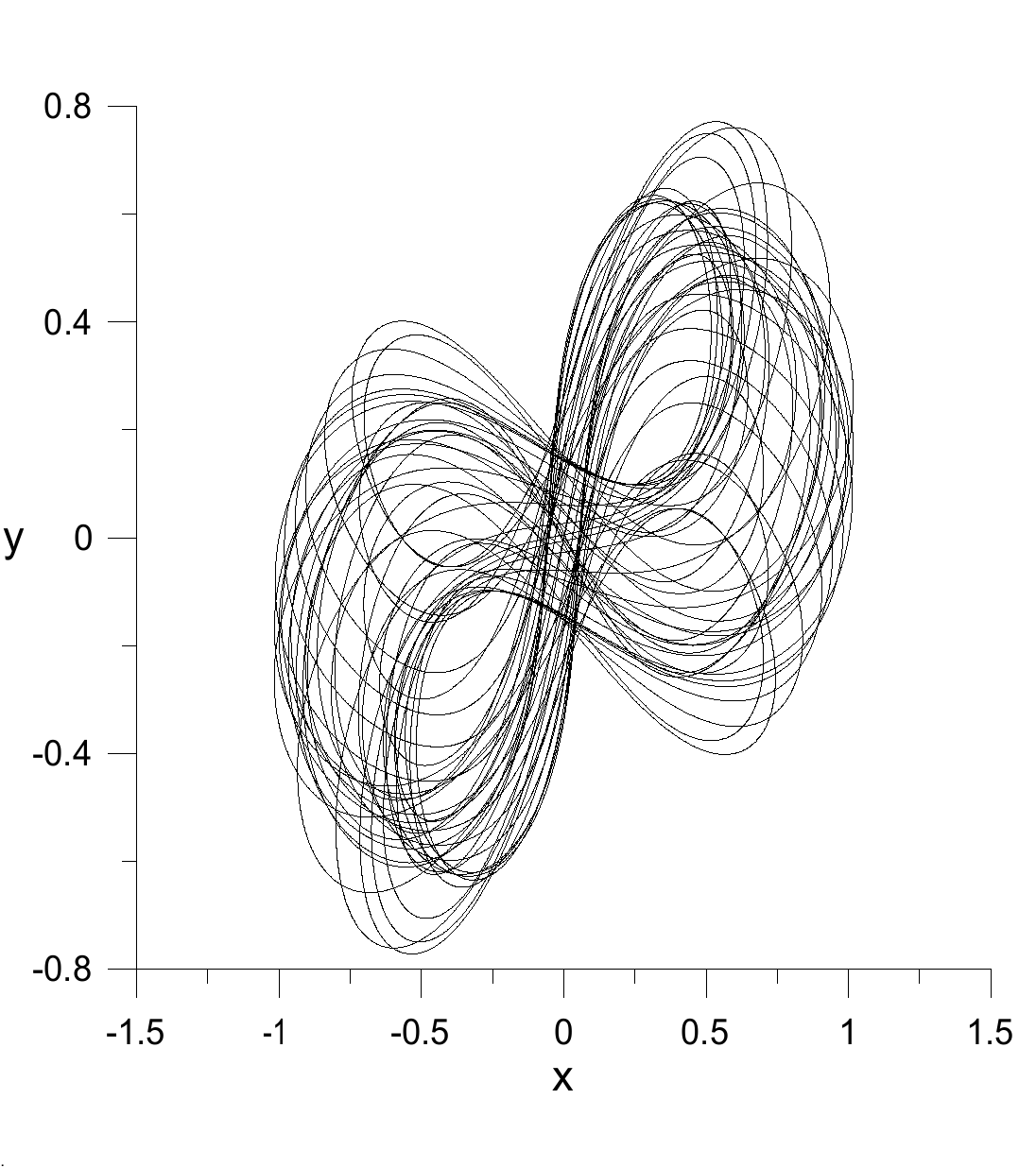}
\caption{Linearly stable choreography (solution 4)}
\label{fig:1}
\end{center}
\end{figure}

\begin{figure}
\begin{center}
\includegraphics[scale=0.52]{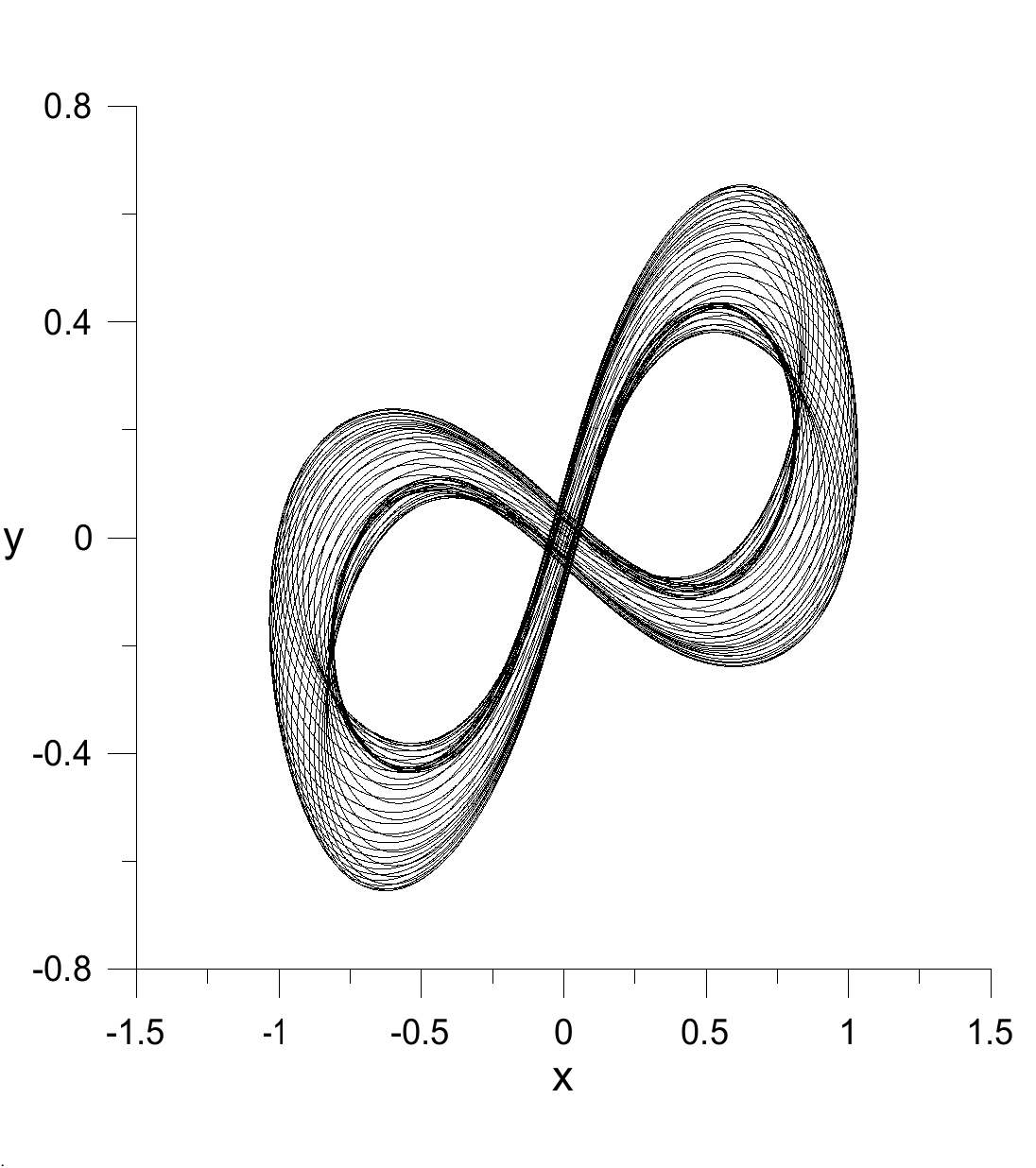}
\caption{Linearly stable choreography (solution 6)}
\label{fig:1}
\includegraphics[scale=0.52]{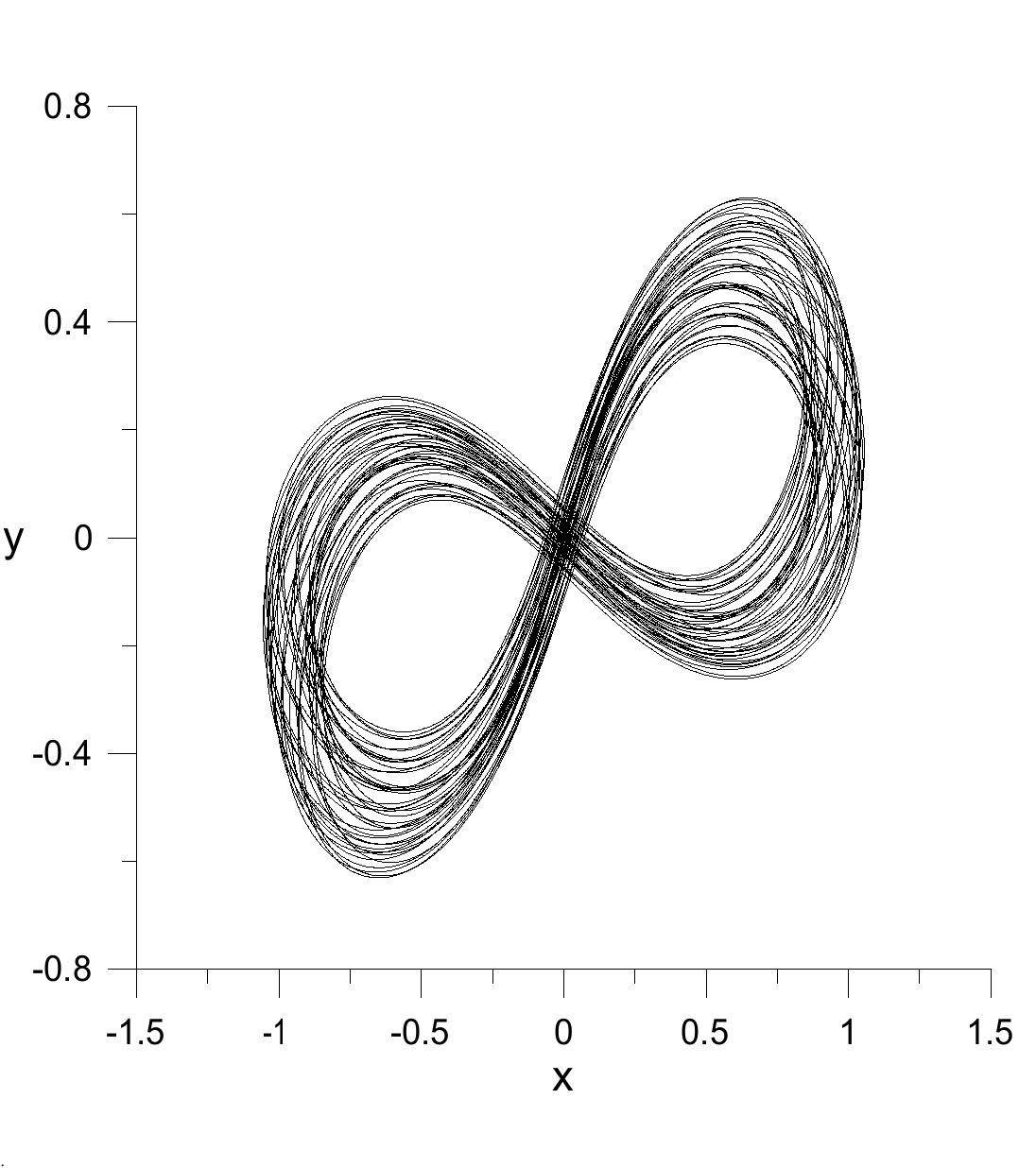}
\caption{Linearly stable choreography (solution 7)}
\label{fig:1}
\end{center}
\end{figure}

\begin{figure}
\begin{center}
\includegraphics[scale=0.52]{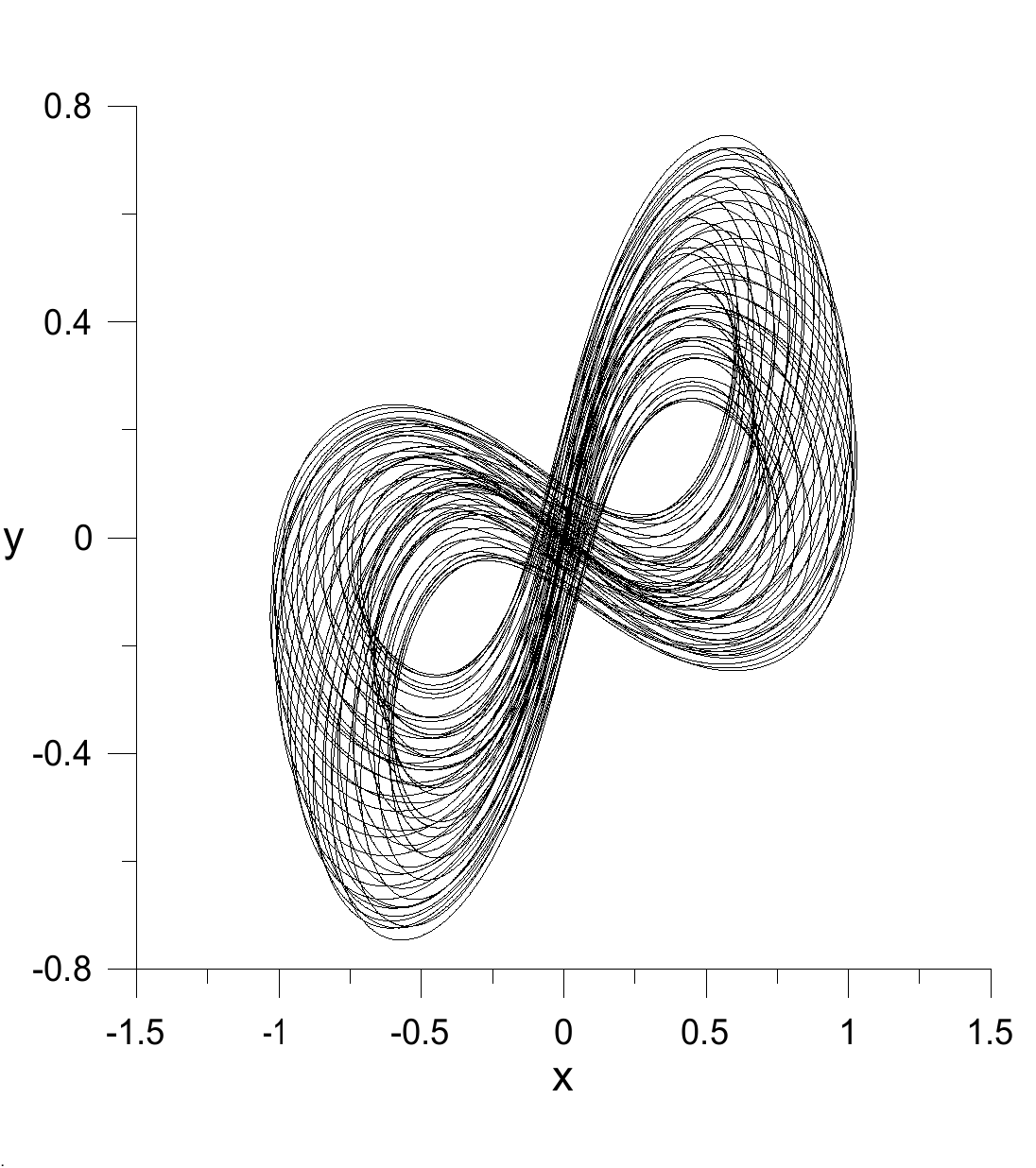}
\caption{Linearly stable choreography (solution 8)}
\label{fig:1}
\includegraphics[scale=0.52]{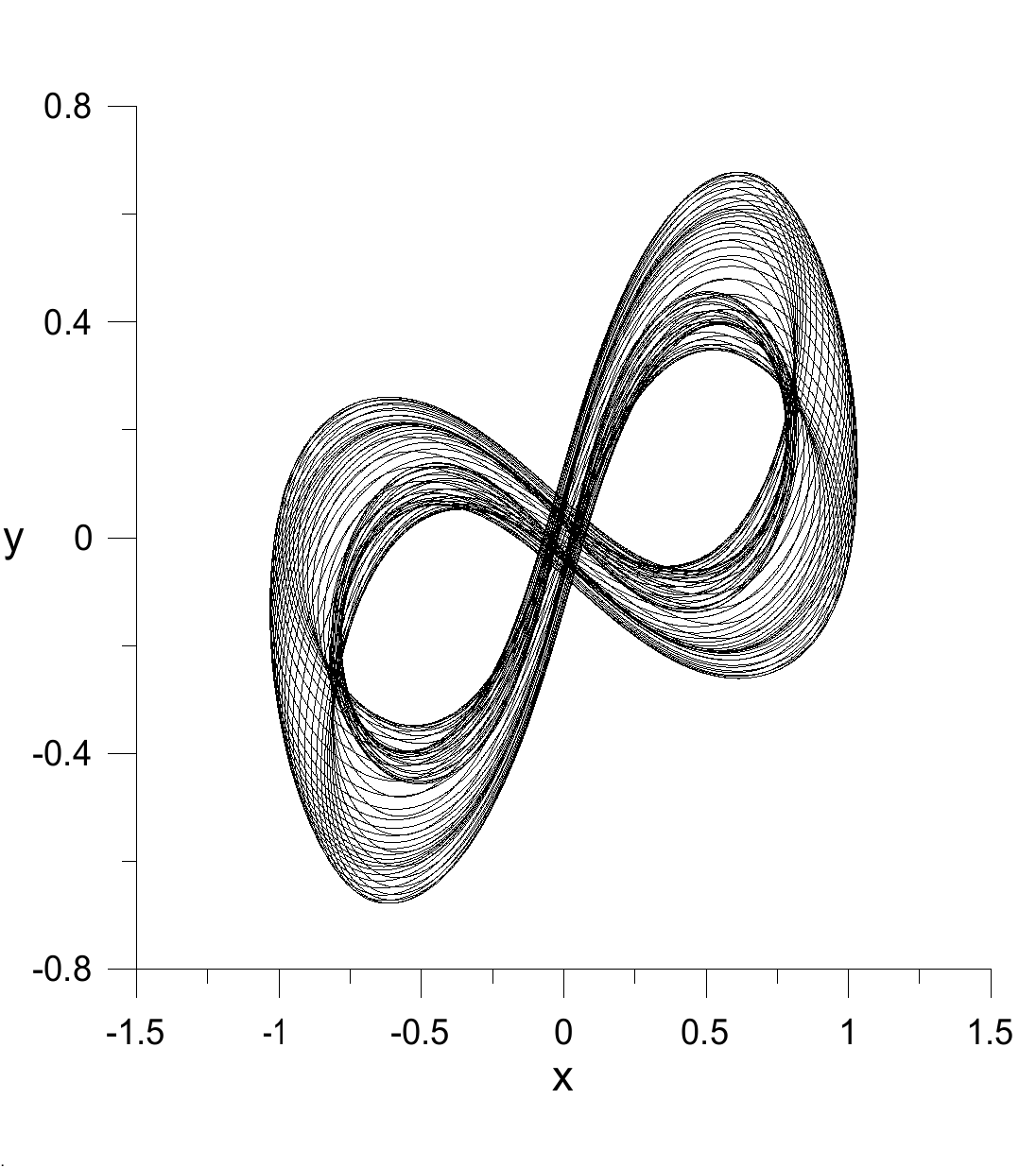}
\caption{Linearly stable choreography (solution 9)}
\label{fig:1}
\end{center}
\end{figure}

\newpage

\end{document}